# A Rolling PID Control Approach and Its Applications

Xiaojun Zhou[1], Chunhua Yang[1], Weihua Gui[1]

1. School of Information Science and Engineering, Central South University, Changsha 410083, China
E-mail: michael.x.zhou@csu.edu.cn

**Abstract:** The canonical proportional-integral-derivative (PID) control approach has been widely used in industrial application due to their simplicity and ease of use. However, its corresponding controller parameters are hard to be adjusted, especially for nonlinear systems. The optimization-based method provides a general framework to find optimal PID controller parameters; nevertheless, several disadvantages exist, for example, it is nontrivial to select an appropriate sample size and it is necessary to obtain the global optimal solution but the optimization problem is non-convex, making it hard to achieve. To alleviate the aforementioned limitations, a rolling PID control approach is proposed in this study, in which, at each rolling period, the PID controller parameters are updated using observable data, which can be classified to data-driven control method. The effectiveness of the proposed approach has been validated by experiments.

## 1 INTRODUCTION

The proportional-integral-derivative (PID) controller and its variants (PI controller for instance) have been widely used in control engineering, and it is reported that over 90% of control loops in industrial applications are of PID-type [1-4]. The advantages of PID controller can be summarized as follows: i) simple control structure but powerful functionality in both transient and steady-state responses; ii) clear physical meaning of controller parameters; iii) easy to be implemented. The main difficulty in PID controller design is to determine three controller parameters, *i.e.*, the proportional gain $K_p$, the integral gain $K_i$ and the derivative gain $K_d$. For single-input single-output (SISO) system, Ziegler-Nichols (ZN) method and its modifications are probably the most known and widely used ones, which are heuristic tuning methods. Another popular approach with similar emphasis for SISO system is the gain and phase margin method. In recent years, optimization-based method for tuning of PID controller parameters has received considerable attention due to its wide generality and flexibility [3-10]. By different specifications, like load disturbance response, set point response, robustness with respect to model uncertainties, the PID controller design problem can be converted into an optimization problem either with constraints or multiple objectives [11-14].

Although the optimization-based PID controllers design method has many benefits, there still exist several limitations in the canonical PID control approach: i) it is nontrivial to select an appropriate sample size $N$; ii) it is necessary to obtain the global optimal solution but the optimization problem non-convex, making it hard to achieve; and iii) it is not robust when model uncertainty exists. To overcome these limitations, a rolling PID control approach is proposed in this study, in which, at each rolling period, the PID controller parameters are updated using observable data, which can be classified to data-driven control method. The feasibility and effectiveness of the proposed rolling PID approach has been demonstrated by several case studies.

The remainder of this paper is organized as follows. In Section 2, a brief review of canonical PID control approach is described, and by discussing its disadvantages, then, a rolling PID control approach is proposed to alleviate these limitations. In Section 3, simulation results are given to illustrate the effectiveness of the proposed approach. Conclusions and future directions are derived in Section 4.

## 2 ROLLING PID CONTROL

### 2.1 Problem statement

In this study, consider the following multivariable SISO discrete-time nonlinear system:

$$\begin{aligned} x(k+1) &= f(k, x(k), u(k), \theta_1) \\ y(k) &= h(k, x(k), u(k), \theta_2) \end{aligned} \quad (1)$$

where $x(k) \in \Re^n$ is the state vector, $u(k) \in \Re$ is the control and $y(k) \in \Re$ is the output vector. $f(\bullet), h(\bullet)$ are given nonlinear functions, and $\theta_1, \theta_2$ are known system parameters. The initial state vector is $x(0) = x_0$.

The goal of the dynamical system is to drive the output vector $y(k)$ to the desired reference signal $y_r(k)$ while satisfying the specified state and control constraints by designing appropriate control approach.

As stated in [8], a typical control approach is composed of two steps: selection of controller structure and tuning of controller parameters. For linear system, linear state or output feedback control law is commonly used. While for system with time delays, uncertainties, external disturbances,

This work is supported by National Nature Science Foundation under (Grant No. 61503416, 61533020, 61533021,61590921)

etc and even nonlinear systems, other types of controller structures are more popular, among which, the PID-type controller is considered as the most popular one in real-world industrial applications. Therefore, the PID-type controller is focused in this study.

**2.2 Canonical PID control approach**

The discrete-time canonical PID control law can be given in *place type* as follows:

$$u(k) = K_p e(k) + K_i \sum_{j=0}^{k} e(j) + K_d \left[ e(k) - e(k-1) \right] \quad (2)$$

where $K_p, K_i, K_d$ are proportional gain, integral gain and derivative gain, respectively; the error $e(k)$ is

$$e(k) = y_r(k) - y(k) \quad (3)$$

here $y_r(k)$ denotes the reference signal and $y(k)$ is the output.

The above PID control law can also be described in the following *increment type*:

$$\begin{aligned} u(k) = u(k-1) &+ K_p \left[ e(k) - e(k-1) \right] \\ &+ K_i e(k) + K_d \left[ e(k) - 2e(k-1) + e(k-2) \right] \end{aligned} \quad (4)$$

To determine the three PID controller parameter vector $K = [K_p, K_i, K_d]$, we define the following optimization problem

$$\begin{aligned} \min_{K} \quad & J(e(k)) \\ \text{s.t.} \quad & g(k, x(k), u(k)) \leq 0 \end{aligned} \quad (5)$$

where $g(k, x(k), u(k)) \leq 0$ is the constraints imposed on the states and control, and $J(e(k))$ is the objective function, which can be expressed as follows under different objective criteria:

$$\begin{cases} J_1(e(k)) = \sum_{k=1}^{N} \left[ e(k) \right]^2 \\ J_2(e(k)) = \sum_{k=1}^{N} k \left[ e(k) \right]^2 \\ J_3(e(k)) = \sum_{k=1}^{N} |e(k)| \\ J_4(e(k)) = \sum_{k=1}^{N} k |e(k)| \end{cases} \quad (6)$$

here $J_1, J_2, J_3, J_4$ are called integral squared error (ISE) criterion, integral time squared error (ITSE) criterion, integral absolute error (IAE) criterion, and integral time absolute error (ITAE) criterion, respectively. Without loss of generality, the ISE criterion is chosen in this study. $N$ is the sample size. Since the interval of integration should be $[0, \infty)$ in continuous time space, $N$ should approach infinity in an ideal situation, which remains as a significant issue for numerical computation.

The schematic diagram of optimization based canonical PID control approach is shown in Fig. 1, and we can find that the PID controller parameters are obtained by solving the optimization problem (5) offline.

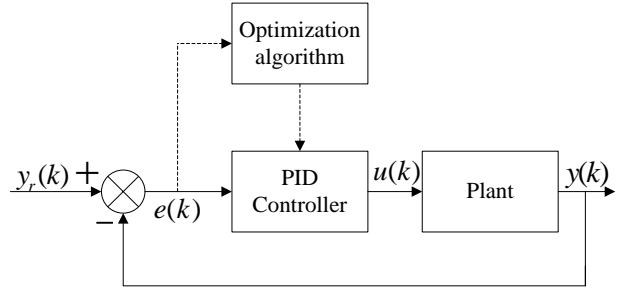

Fig 1. Schematic diagram of optimization based canonical PID control approach.

**2.3 Rolling PID control approach**

In the optimization-based canonical PID control approach, the PID controller parameters are determined offline, and they keep constant in the implementation stage. There exist several disadvantages: i) it is nontrivial to select an appropriate sample size *N*; ii) it is necessary to obtain the global optimal solution but the optimization problem is non-convex, making it hard to achieve; iii) it is not robust when model uncertainty exists.

To eliminate the limitations, an optimization based rolling PID control approach is proposed in this paper, in which, at each rolling period, the PID controller parameter vector is updated using observable data. As shown in Fig.2, the schematic diagram of optimization based rolling PID control approach is composed of two stages: model updating stage (I), and PID controller parameters updating stage (II). To be more specific, at the model updating stage (I), an approximate model is established using available historical measurements, while at the controller parameter vector updating stage (II), the PID controller parameters are updated by solving a rolling optimization problem, which resembles the finite horizon optimal control problem. In a certain sense, the proposed rolling PID approach also has some perspective on the future.

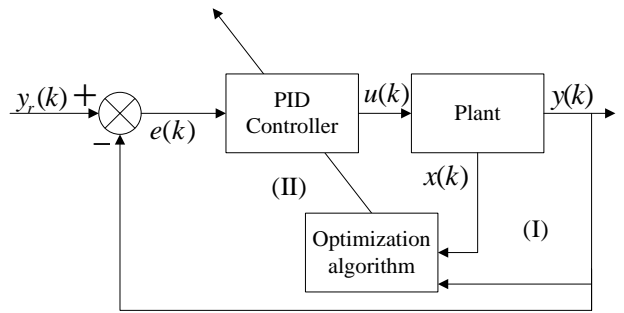

Fig 2. Schematic diagram of optimization based rolling PID control approach.

When using the rolling PID control approach, the objective criteria should be slightly modified as

$$\begin{cases} J_1(e(k)) = \sum_{j=k}^{k+N-1} [e(j)]^2 \\ J_2(e(k)) = \sum_{j=k}^{k+N-1} (j-k+1)[e(j)]^2 \\ J_3(e(k)) = \sum_{j=k}^{k+N-1} |e(j)| \\ J_4(e(k)) = \sum_{j=k}^{k+N-1} (j-k+1)|e(j)| \end{cases} \quad (7)$$

where $N$ can be regarded as the predictive horizon length. The detailed steps of the proposed optimization based rolling PID control approach are given as follows:

**Step 0:** $s=0$, select predictive horizon length $N$ and sample size $M$ by experience; by solving the optimization problem (5) using an optimization method (either local or global) with initial state vector $x_0$, the first PID controller parameter vector $K(s)$ will be obtained.

**Step 1**: $s=s+1$, keep the vector $K(s-1)$ constant until $M$ sampling periods, update the model using data $x(k)$, $x(k-1)$, …, $x(k-M+1)$, $y(k)$, $y(k-1)$, …, $y(k-M+1)$.

**Step 2**: by solving the optimization problem (5) using an optimization method (either local or global) with initial state vector $x(k)$ and sample size $N$, the updated PID controller parameter vector $K(s)$ will be obtained.

**Step 3**: repeat Step 1 until some termination criteria are met.

**Remark 1**: The proposed rolling PID (rPID) control approach resembles the model predictive control (MPC) [16] but distinguishes from it distinctly in the control law. In MPC, the control sequence is considered as independent decision variable which is obtained by solving a finite horizon optimal control problem, and only the first control vector in the control sequence is applied to the plant. While in rPID, the decision variable is PID control parameter vector, which is obtained by solving a series of non-convex optimization problems, and the control law is generated by the PID type control.

## 3 SIMULATION RESULTS

In this section, several examples are illustrated to show the effectiveness of the proposed rolling PID approach.

Let us first consider the following discrete-time nonlinear system as a motivating example:

$$x_1(k+1) = \theta_1 x_1(k) x_2(k), \; x_1(0) = 1,$$
$$x_2(k+1) = \theta_2 x_1^2(k) + u(k), \; x_2(0) = 1,$$
$$y(k) = \theta_3 x_2(k) - \theta_4 x_1^2(k),$$

where $\theta_1, \cdots \theta_4$ are system parameters. The goal is to drive the system output $y$ to the desired reference signal $y_r$.

**Case 1**: $y_r=2$, $[\theta_1,\theta_2,\theta_3,\theta_4]=[0.5,0.3,1.8,0.9]$ without model updating

In this case, without mode updating means that the nonlinear system can be accurately obtained. Assuming that the lower and upper bounds of PID controller parameters are

$$0 \leq K_p, K_i, K_d \leq 10.$$

Let us set $N = 10$, $M = 10$, the initial control parameter vector $K_0$=[0.1,0.1,0.1]. By using a local optimization solver based on sequential quadratic programming (SQP) method, the first PID controller parameter vector can be obtained as $K(1) = [0.0707, 0.3634, 0.1498]$, and then $[x_1(10), x_2(10)] = [0.0005, 1.0971]$. Sequentially, the corresponding results can be obtained as shown in Table 1 and Fig 3. It can be found that the controller parameter vector, states and output are kept constant when $s = 3$, and the desired reference signal is achieved since then.

Table 1. Results obtained by the proposed rPID control without model updating ($N = 10$, $M = 10$)

| $s$ | $K(s)$ | $[x_1(M*s), x_2(M*s)]$ | $y(M*s)$ |
|---|---|---|---|
| 1 | [0.0707  0.3634  0.1498] | [0.0005  1.0971] | 1.8785 |
| 2 | [0.0635  0.2333  0.0640] | [0.0000  1.1110] | 1.9990 |
| 3 | [0.0635  0.2333  0.0640] | [0.0000  1.1111] | 2.0000 |
| 4 | [0.0635  0.2333  0.0640] | [0.0000  1.1111] | 2.0000 |

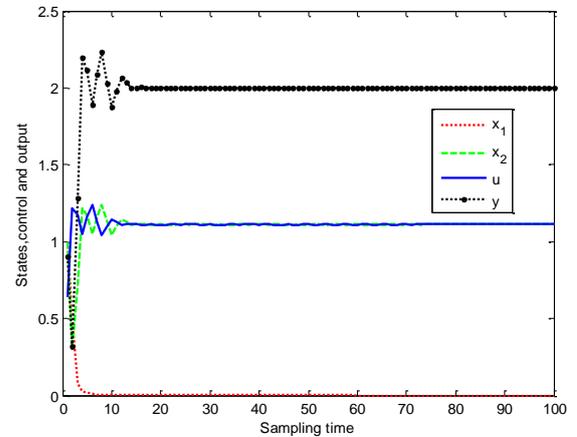

Fig 3. Trajectories of states, control and output ($N = 10$, $M = 10$, Case 1).

Next, another pair $(N,M) = (30,30)$ is used. In this situation, the first PID controller parameter vector can be obtained as $K(1) = [0.0697, 0.3625, 0.1472]$, and then $[x_1(10), x_2(10)] = [0.0000, 1.1111]$. Subsequent results are listed in Table 2 and Fig. 4. Compared with the first pair when $(N,M) = (10,10)$, fewer least rolling times are needed ($s=2$); however, its trajectory of output tends to steady much late.

Table 2. Results obtained by the proposed rPID control without model updating ($N = 30$, $M = 30$)

| $s$ | $K(s)$ | $[x_1(M*s), x_2(M*s)]$ | $y(M*s)$ |
|---|---|---|---|
| 1 | [0.0697  0.3625  0.1472] | [0.0000  1.1111] | 1.9986 |
| 2 | [0.0342  0.3288  0.1221] | [0.0000  1.1111] | 2.0000 |
| 3 | [0.0342  0.3288  0.1221] | [0.0000  1.1111] | 2.0000 |

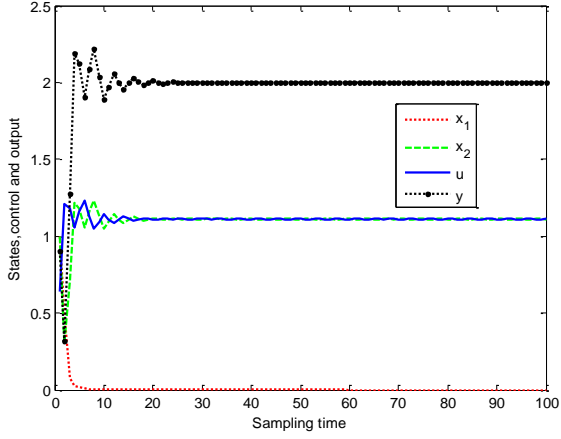

Fig 4. Trajectories of states, control and output ($N = 30$, $M = 30$, Case 1).

**Case 2**: $y_r=2$, $[\theta_1,\theta_2,\theta_3,\theta_4]=[0.5,0.3,1.8,0.9]$ with model updating

In this case, the following discrete-time linear system is used to approximate the nonlinear model.

$$\hat{x}(\hat{k}+1) = A\hat{x}(\hat{k}) + B\hat{u}(\hat{k}), \hat{k} = k, k-1, \cdots k-M+1$$
$$\hat{y}(\hat{k}) = C\hat{x}(\hat{k}),$$

where $A = \begin{bmatrix} a_{11} & a_{12} \\ a_{21} & a_{22} \end{bmatrix}$, $B = \begin{bmatrix} 0 \\ b_1 \end{bmatrix}$, $C = \begin{bmatrix} c_1 & c_2 \end{bmatrix}$ are unknown parameters needed to be determined utilizing $M$ sampling data at each rolling period. By using least squares method, the following optimization problem is formulated

$$\arg\min_{A,B,C} \sum_{\hat{k}=k-M+1}^{k} \left\| \hat{x}(\hat{k}+1) - A\hat{x}(\hat{k}) - B\hat{u}(\hat{k}) \right\|_2^2 + \left| \hat{y}(\hat{k}) - C\hat{x}(\hat{k}) \right|^2$$

Let us set $N = 10$, $M = 10$, the initial control parameter vector $K_0=[0.1,0.1,0.1]$. By solving the above convex optimization problem, parameters in the approximate linear model can be obtained as follows

$$A = \begin{bmatrix} 0.4236 & 0.0056 \\ 0.1588 & 0.1392 \end{bmatrix}, B = \begin{bmatrix} 0 \\ 0.8711 \end{bmatrix},$$
$$C = \begin{bmatrix} -0.8041 & 1.7975 \end{bmatrix}$$

Next, the PID controller parameter vector is updated based on the approximate linear model by solving the optimization problem (7). At this stage, $K(1) = [0.0944,0.2340,0.0930]$. Sequentially, the corresponding results can be obtained as shown in Table 3 and Fig 5. It can be found that the controller parameter vector, states and output are kept constant when $s = 4$, and the desired reference signal is achieved since then. Compared with the case when no model updating is considered, it is observed that no overshoot is observed, and the trajectories are much steadier because no oscillation occurs.

**Remark 2**: The meaning of the initial control parameter vector here is different from that of Case 1, in which, it is used as starting point for the SQP method. While in this case, it is implemented as the first PID controller parameter vector.

Table 3. Results obtained by the proposed rPID control with mode updating ($N = 10$, $M = 10$)

| s | K(s) | | | [$x_1$(M*s), $x_2$(M*s)] | | y(M*s) |
|---|---|---|---|---|---|---|
| 1 | [0.1000 | 0.1000 | 0.1000] | [0.0000 | 0.8916] | 1.6049 |
| 2 | [0.0944 | 0.2340 | 0.0930] | [0.0000 | 1.1092] | 1.9966 |
| 3 | [0.0627 | 0.2308 | 0.0624] | [0.0000 | 1.1110] | 1.9999 |
| 4 | [0.0627 | 0.2308 | 0.0624] | [0.0000 | 1.1111] | 2.0000 |
| 5 | [0.0627 | 0.2308 | 0.0624] | [0.0000 | 1.1111] | 2.0000 |

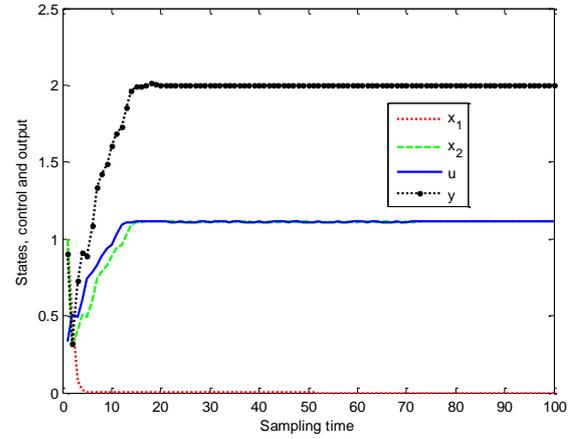

Fig 5. Trajectories of states, control and output ($N = 10$, $M = 10$, Case 2).

Now, let us set $N = 30$, $M = 30$, similar results can be found in Table 4 and Fig 6 but with more stable trajectories.

Table 4. Results obtained by the proposed rPID control with model updating ($N = 30$, $M = 30$)

| s | K(s) | | | [$x_1$(M*s), $x_2$(M*s)] | | y(M*s) |
|---|---|---|---|---|---|---|
| 1 | [0.1000 | 0.1000 | 0.1000] | [0.0000 | 1.1071] | 1.9928 |
| 2 | [0.0458 | 0.1110 | 0.0231] | [0.0000 | 1.1111] | 2.0000 |
| 3 | [0.0512 | 0.1310 | 0.0241] | [0.0000 | 1.1111] | 2.0000 |
| 4 | [0.0512 | 0.1310 | 0.0241] | [0.0000 | 1.1111] | 2.0000 |

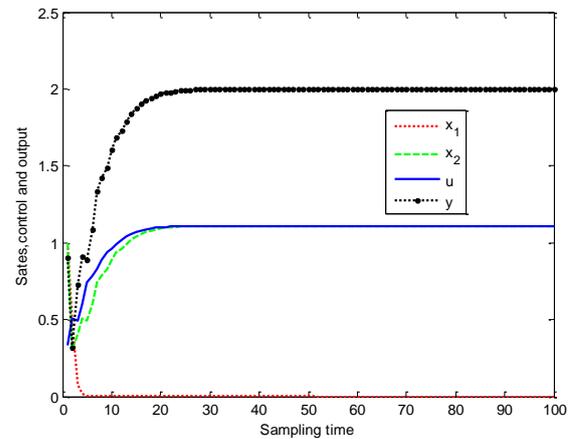

Fig 6. Trajectories of states, control and output ($N = 30$, $M = 30$, Case 2).

## 4 CONCLUSION AND FUTURE DIRECTIONS

A rolling PID control approach is proposed in this study, the core of which is to update the PID controller parameters periodically using available measured data. To be more specific, an approximate model is established using historical measurements, and then the PID controller parameters are updated by solving a rolling optimization problem. The effectiveness of the proposed approach is verified by simulation results.

On the other hand, for control engineering problems, states and control constraints are commonly met in real-world applications. A future direction will include the extension of the rolling PID approach to constrained nonlinear system. Furthermore, only the tracking problem is considered in this study, other control problem will be studied. In the meanwhile, the philosophy (theory) behind the proposed methodology should be well addressed.